\documentclass[12pt]{article}
\usepackage{fullpage}
\usepackage{amsfonts}
\usepackage{amsmath}
\ifx\pdfoutput\undefined \csname newcount\endcsname\pdfoutput \fi
\ifcase\pdfoutput
  \usepackage[dvips]{graphicx}
  \DeclareGraphicsExtensions{.eps}
\else
  \usepackage[pdftex]{graphicx}
  \DeclareGraphicsExtensions{.pdf,.jpg}
\fi

\def\R{R}
\def\qed{{$\Box$}}

\newcommand{\flushpar}{\vspace{1em}\noindent}

\newtheorem{Proposition}{Proposition}
\newtheorem{Lemma}{Lemma}
\newtheorem{Corollary}{Corollary}
\newtheorem{Problem}{Problem}
\newenvironment{Proof}{{\noindent\em Proof:\quad}}{\mbox{}\qed\\[1em]\mbox{}}
\newenvironment{Remark}{{\noindent\em Remark:\quad}}
{\mbox{}$\Diamond$\\[1em]\mbox{}}
\newenvironment{Example}{{\noindent\em Example:\quad}}
{\mbox{}$\Diamond$\\[1em]\mbox{}}

\newcommand{\har}[2]{{!({#1},{#2})}}

\begin{document}

\vskip 1in
\bigbreak \centerline {\bf Zero-Preserving Iso-spectral
Flows Based on Parallel Sums}
\par\bigbreak

\centerline{by}
\par\medbreak
\centerline{Kenneth R. Driessel}
\centerline{Mathematics Department}
\centerline{Colorado State University}
\centerline{Fort Collins, Colorado, USA}
\centerline{email: driessel\makeatletter{@}math.colostate.edu}
\par\medbreak
\centerline{and}
\par\medbreak
\centerline{Alf Gerisch}
\centerline{Fachbereich Mathematik und Informatik}
\centerline{Martin-Luther-Universit\"at Halle-Wittenberg}
\centerline{06099 Halle (Saale), Germany}
\centerline{email: gerisch\makeatletter{@}mathematik.uni-halle.de}
\par\medbreak

\bigbreak
\centerline{October, 2001, Revised in 2004 and 2005}

\vskip 1in
\bigbreak


\newpage

\subsection*{Abstract}
Driessel~[\emph{Computing canonical forms using flows}, Linear Algebra and
Its Applications 2004] introduced the notion of quasi-projection onto
the range of a linear transformation from one inner product space into
another inner product space. Here we introduce the notion of
quasi-projection onto the intersection of the ranges of two linear
transformations from two inner product spaces into a third inner
product space. As an application, we design a new family of
iso-spectral flows on the space of symmetric matrices that preserves
zero patterns. 
We discuss the equilibrium points of these flows.
We conjecture that
these flows generically converge to diagonal matrices.  We perform
some numerical experiments with these flows which support this
conjecture. We also compare our zero preserving flows with the Toda
flow.

\vspace{-5mm}

\tableofcontents

\newpage

\section{Introduction.}

Let $\Delta$ be a set of pairs $(i,j)$ of integers between 1 and $n$
which satisfies the following conditions: (1) for $i=1,2,...,n$, the
diagonal pair $(i,i)$ is in $\Delta$, and (2) if the pair $(i,j)$ is
in $\Delta$ then so is the symmetric pair $(j,i)$. We regard
$\Delta$ as a (symmetric) sparsity {\bf pattern of interest}
of nonzero entries for
matrices. In particular, let $Sym(n)$ denote the vector space of
symmetric, $n\times n$, real matrices and let $Sym(\Delta)$ denote the
subspace of $Sym(n)$ consisting of the symmetric matrices which are
zero outside the pattern $\Delta$; in symbols,
$$
Sym(\Delta) :=
\{ X \in Sym(n) : X(i,j) \neq 0 \text{ implies } (i,j) \in \Delta \} .
$$

In this report we consider the following task: Find flows in the space
$Sym(\Delta)$ which preserve eigenvalues and converge to diagonal
matrices. We can describe this task more precisely as follows: With an
$n \times n$ symmetric matrix $A$, we associate the {\bf iso-spectral
surface}, $Iso(A)$, of all symmetric matrices which have the same
eigenvalues as $A$. By the spectral theorem, we have
$$
Iso(A) := \{QAQ^T : Q \in O(n) \}
$$
where $O(n)$ denotes the group of orthogonal matrices.

We shall use the Frobenius inner product on matrices; recall that
it is defined by
$\langle X,Y \rangle:= Trace(X Y^T)$.  With a symmetric matrix $D$, we
associate a real-valued `objective' function
$$
f:= Sym(n) \to \R: X\mapsto (1/2)\langle X-D,X-D \rangle .
$$
Note that $f$ is a measure of the distance from $X$ to $D$.  We shall
consider the following constrained optimization problem:

\begin{Problem}
Given $A\in Sym(\Delta)$, minimize $f(X)$ subject to
the constraints $X \in Iso(A)$ and $X \in Sym(\Delta)$.
\end{Problem}
In particular, we shall describe a flow on the surface $Iso(A)
\cap Sym(\Delta)$ which solves this problem in the sense that it
usually converges to a local minimum.

Here is a summary of the contents of this report.

In the next section which is entitled ``Quasi-projection onto the
intersection of two subspaces'', we present some theoretical
background material. Driessel[2004] introduced the notion of
quasi-projection onto the range of a linear transformation from
one inner product space to another. In this section we introduce
the notion of quasi-projection onto the intersection of the ranges
of two linear transformations $A$ and $B$ from two inner product 
spaces into a third inner product space.
We use the notation $!(A,B)$ to denote our quasi-projection operator.
We show that $!(A,B)=2A(A+B)^+B$ where the superscript $+$ denotes the
Moore-Penrose pseudo inverse operation.

\begin{Remark}
If $A$ and $B$ are invertible then
$$
\har{A}{B} = 2A(A+B)^{-1}B = 2(A^{-1}+B^{-1})^{-1} .
$$
This operator is called the ``harmonic mean'' of the operators $A$ and $B$.
See, for example, Kubo and Ando~[1980].
They use the the infix notation $A!B$ to denote
the harmonic mean of $A$ and $B$ where $A$ and $B$ are positive operators
on a Hilbert space.
After we wrote this paper in 2001, Chandler Davis told us about this
paper by Kubo and Ando. This paper led us to the following papers:
Anderson and Duffin[1969], Anderson[1971], Anderson and
Schreiber[1972], Anderson and Trapp[1975]. In particular, Anderson
and Duffin define the ``parallel sum'' of semi-definite matrices
$A$ and $B$ by the formula $A(A+B)^+B$ and denote it by $A:B$.  
We discovered that most of the results in Section~2 appear 
scattered in these earlier papers (but usually with different
proofs). In order to keep this paper somewhat self-contained we
retained our proofs. 
\end{Remark}

In the third section which is entitled ``An iso-spectral flow
which preserves zeros'', we describe an application of the
quasi-projection method. In particular, we describe how we used
this method to design a new flow corresponding to the optimization
problem described above. We conjecture that this flow generically
converges to a symmetric matrix $E$ that commutes with $D$. Note
that if we choose $D$ to be a diagonal matrix with distinct
diagonal entries then $E$ commutes with $D$ iff $E$ is a diagonal
matrix.
(For background material on differential equations see, for example,
Hirsch, Smale and Devaney~[2004].)

In the fourth section which is entitled ``Numerical results'', we
describe our implementation of our iso-spectral zero-preserving
flow in Matlab. We also describe several numerical experiments
that we performed using this computer program. 

In all our experiments this flow converges (sometimes slowly) to a 
diagonal matrix.
Consequently these experiments provide evidence for the conjecture
described above. We do not claim our program to be competitive
with standard methods used to compute eigenvalues. But we hope our
ideas will lead eventually to practical, competitive methods for
finding eigenvalues of some classes of structured matrices.

In a first appendix which is entitled ``Comparing projections and
quasi-projections'', we describe the origin of the quasi-projection
method. In particular, we review a standard method of projection
onto the intersection of the ranges of two linear maps. We show
how quasi-projection arises by simplifying this standard
projection procedure. We also argue that quasi-projection is
simpler, more direct and more robust than projection.

In a second appendix which is entitled ``On the Toda flow'', we
indicate our current geometrical view of the so-called Toda flow
or QR flow. Most of the results in this appendix are known. We
present these results to show the analogies between the
iso-spectral Toda flow and our iso-spectral, zero-preserving flow.
These analogies provided the basis for our development of these
new flows.
(We have repeated some of the definitions of
notation in this appendix. We want to make this appendix
self-contained. We hope that a reader can understand it without
knowledge of the rest of this report.)

For another example of a structured
iso-spectral flow see Fasino~[2001].

\section{Quasi-projection onto the intersection of two subspaces.}

In this section we shall present some theoretical background material
concerning quasi-projections. We shall apply this material in the next
section.  Let $V$ be a finite-dimensional, real inner product
space. We use $\langle x,y \rangle$ to denote the inner product of two
elements of $V$. Let $A:V\to V$ and $B:V\to V$ be (self-adjoint)
positive semi-definite linear operators on $V$. For any vector $c$ in
$V$, consider the following system of linear equations for $u$ and
$\lambda$ in $V$:
\begin{align*}
u - A\lambda &= Ac \tag{ q1} \\
(A+B)\lambda &= (B-A)c \tag{ q2}
\end{align*}
We call these equations the {\bf quasi-projection equations} determined
by $A$, $B$ and $c$.

\begin{Remark}
In this section we usually assume that $A$ and $B$ are two positive
semi-definite operators on a finite dimensional space. These
assumptions simplify the analysis considerably. They will be
obviously satisfied in the application considered below. However,
many of the results in this section are true in more general
settings.
\end{Remark}

Note that (q1) is equivalent to the following condition:
\begin{equation}
u=A(\lambda + c) \tag{eq1} .
\end{equation}
Hence $u$ is in the range of $A$. Also note that (q2) is equivalent to
the following condition:
\begin{equation}
A(\lambda + c) = B(-\lambda +c) \tag{eq2} .
\end{equation}
Hence $u$ is also in the range of $B$. Thus we see that $u$ is in the
intersection of the range of $A$ and the range of $B$.

\begin{Remark}
We sometimes use $f.x$ or $fx$ in place of $f(x)$ to indicate
function application. We do so to reduce the number of
parentheses. We also use association to the left. For example,
$D(\omega.A).I.K$ means evaluate $\omega$ at $A$ to get a
function, differentiate this function, evaluate the result at $I$
to get a linear function, and finally evaluate at $K$.
We adapted this notation from the programming language C
(in which such a dot notation is used in connection with data structures).
\end{Remark}

We shall use the following lemma repeatedly.

\begin{Lemma}
If $A$ and $B$ are positive semi-definite
operators then
\begin{align*}
Kernel(A+B) &= Kernel.A \cap Kernel.B \,,\\
Range(A+B) &= Range.A + Range.B .
\end{align*}
\end{Lemma}

\begin{Proof}
 If $Az = Bz = 0$ then $(A+B)z=0$. Now assume $(A+B)z=0$.
Then $0=\langle z,(A+B)z \rangle = \langle z, Az \rangle + \langle z,
Bz \rangle$.  Since $A$ and $B$ are positive semi-definite, we get
$0=\langle z, Az \rangle = \langle z, Bz \rangle$ and hence
$0=Az=Bz$. The second equation of this lemma is obtained from the
first one by taking orthogonal complements.
\end{Proof}

The following proposition shows that the vector $u$ is uniquely
determined by the quasi-projection equations.

\begin{Proposition}[Uniqueness] Let $A$ and $B$ be positive
semi-definite operators. For any $c\in V$, if $(u_1,\lambda_1)$ and
$(u_2,\lambda_2)$ are solutions of the quasi-projection equations (q1)
and (q2) then $u_1=u_2, A\lambda_1 = A\lambda_2$ and $B\lambda_1 =
B\lambda_2$.
\end{Proposition}
\begin{Proof}
Let $u:=u_1-u_2$ and $\lambda:=\lambda_1 -\lambda_2$. Then we have
$u-A\lambda = 0$ and $(A+B)\lambda=0$. By Lemma~1 we get $A\lambda
= B\lambda = 0$. Then $u = A\lambda = 0$.
\end{Proof}

The following proposition shows that solutions of the quasi-projection
equations always exist.

\begin{Proposition}[Existence]
Let $A$ and $B$ be positive
semi-definite operators. For all $c\in V$, there exist $u$ and
$\lambda$ in $V$ satisfying the quasi-projection equations (q1) and
(q2).
\end{Proposition}

\begin{Proof}
It clearly suffices to show that there is a $\lambda$ in
$V$ such that $(A+B)\lambda = (B-A)c$. In other words, we need
to see that $(B-A)c \in Range(A+B) = Range.A + Range.B$. For this we
simply note $(B-A)c = A(-c) + Bc \in Range.A + Range.B$.
\end{Proof}

Let $\har{A}{B}: V \to V$ denote the linear operator on $V$ which
maps a vector $c$ to the unique vector $u$ which satisfies the
following condition: There exists $\lambda \in V$, such that the
pair $(u,\lambda)$ satisfies the quasi-projection equations (q1)
and (q2). We call the vector $u=\har{A}{B}.c$ the {\bf
quasi-projection} of $c$ onto the intersection of $Range.A$ and
$Range.B$.
Following Anderson and Duffin~[1969] we call $!(A,B)$ the 
{\bf parallel sum} of $A$ and $B$ (even though there is a difference of a
factor of $2$).

For any linear map $M$ between inner product spaces let $M^*$ denote
the {\bf adjoint} map which is defined by the following condition: for
all $x$ in the domain of $M$ and all $y$ in the codomain of $M$,
$\langle Mx,y \rangle = \langle x,M^*y \rangle$.  (Halmos~[1958] uses
this notation for the adjoint.)  The following proposition shows how
quasi-projection behaves with respect to congruence.

\begin{Proposition}[Congruence]
Let $M:V\to V$ be any invertible
linear map. Then $M (\har{A}{B}) M^* = \har{MAM^*}{MBM^*} $.
\end{Proposition}
\begin{Proof}
The pair of equations (eq1) and (eq2) is equivalent to the following
pair:
\begin{align*}
Mu&=MAM^*(M^*)^{-1}(c+\lambda), \\
MAM^*(M^*)^{-1}(c+\lambda) &= MBM^*(M^*)^{-1}(c-\lambda) .
\end{align*}
Hence, for all $c$ in $V$, we have
$$
M(!(A,B))c = !(MAM^*,MBM^*) (M^*)^{-1} c .
$$
\end{Proof}

Let $U$ and $V$ be inner product spaces and let $L : U \to V$ be a linear
map.
We use $L^+$ to denote the Moore-Penrose pseudo-inverse of $L$. (See,
for example, Lawson and Hanson~[1974].) We list the following
properties of the pseudo-inverse
$$
L^{*+} = L^{+*},\quad L L^+ L = L,\quad L^+ L L^+ = L^+ ,
$$
and note that $L L^+$ is the projection of $V$ onto $Range.L$ and
$L^+ L$ is the projection of $U$ onto $Range.L^*$.

\begin{Lemma}
Let $A$ and $B$ be positive semi-definite operators on an inner 
product space $V$. Then
$$
A = A(A+B)(A+B)^+ = A(A+B)^+(A+B) = (A+B)(A+B)^+A = (A+B)^+(A+B)A .
$$
\end{Lemma}
\begin{Proof}
Note that $P:=(A+B)(A+B)^+ = (A+B)^+(A+B)$ is the projection of $V$
onto the range of $A+B$. In particular, by Lemma 1, for all $x$ in the
range of A, we have $Px=x$. Also note that $V=Range.A \oplus Kernel.A$
since $(Range.A)^\bot = Kernel.A^* = Kernel.A$.

%
Now consider any $x\in V$. Note that $PAx=Ax$ because $Ax$ is in the 
range of $A$ which is a subset of the range of $A+B$.
Hence $A=PA$. Since $A$ is self-adjoint we also
have $A=AP$.
\end{Proof}

The following proposition is our main result concerning
quasi-projections. We shall use it below to design zero preserving
flows.

\begin{Proposition}[Quasi-Projection Formulas]
Let $A$ and $B$ be positive semi-definite operators. Then the
quasi-projection operator is given by the following formulas:
$$
\har{A}{B}= 2A(A+B)^+B = 2B(A+B)^+A .
$$
Furthermore, the quasi-projection operator is positive semi-definite.
Its range equals the intersection of the range of $A$ and the range of
$B$ and its kernel equals the sum of the kernel of $A$ and the kernel
of $B$; in symbols,
\begin{align*}
Range (\har{A}{B}) &= Range.A\cap Range.B \\
Kernel (\har{A}{B}) &= Kernel.A + Kernel.B .
\end{align*}
\end{Proposition}

\begin{Proof}

\flushpar Claim: $\har{A}{B}=2A(A+B)^+B$ .

We take $\lambda := (A+B)^+(B-A)c$. This $\lambda$ satisfies the
quasi-projection equation (q2). Substituting in equation (q1), we get
$ u=\har{A}{B}c = (A(A+B)^+(B-A) + A)c .$
Using Lemma 2, we get
$ A(A+B)^+(B-A) + A = 2A(A+B)^+B $ .

\flushpar Claim: $A(A+B)^+B = A - A(A+B)^+A$ .

Using Lemma 2 again yields
$
A - A(A+B)^+A = A(A+B)^+(A+B) - A(A+B)^+A = A(A+B)^+B .
$

\flushpar Claim: The map $A(A+B)^+B$ is self-adjoint.

Use the previous claim and the fact that $(A+B)^+=(A+B)^{*+} =
(A+B)^{+*}$.

\flushpar Claim: $A(A+B)^+B = B(A+B)^+A$ .

Use the previous claim and $(A(A+B)^+B)^* = B(A+B)^+A$.

\flushpar Claim: $Kernel(\har{A}{B}) = Kernel.A + Kernel.B$ .

By the formulas for the quasi-projection, we see that its kernel
contains $Kernel.A$ and $Kernel.B$ and hence $Kernel.A+Kernel.B$. We
need to prove the other inclusion; in other words, we want to see that
the following statement is true:

$$
\forall z\in Kernel.(\har{A}{B}), \exists x,y \in V,
z=x+y,Ax=0, By=0 .
$$
Consider any $z$ satisfying  $0=\har{A}{B}z=2A(A+B)^+Bz$. Take
$x:=(A+B)^+Bz$. Note $Ax=0$. Using Lemma 2 again we also have $Bx
= (A+B)x = (A+B)(A+B)^+Bz = Bz$. Hence $B(z-x)=0$. We can take
$y:=z-x$.

\flushpar Claim: $Range(\har{A}{B}) = Range.A \cap Range.B$

Take orthogonal complements of the previous claim.

\flushpar Claim: The map $\har{A}{B}$ is positive semi-definite.

Note that the range of $A+B$ is an invariant subspace of
$\har{A}{B}$. Clearly we only need to see that the restriction of
$\har{A}{B}$ to this range is positive semi-definite. %
%
Consequently we assume that $V=Range(A+B)$. In this case we have
$\har{A}{B}=2A(A+B)^{-1}B=2B(A+B)^{-1}A$. We now view $A$ and $B$ as
matrices. Since $A+B$ is positive definite and $A$ is
self-adjoint, we can simultaneously diagonalize these two matrices
by a congruence. (See, for example, Bellman~[1970].)
In particular, there is an invertible
matrix $M$ and a diagonal matrix $D:=diag(a_1^2,\dots,a_n^2)$ such
that $M(A+B)M^*=I$ and $MAM^*=D$.  We see from these equations
that $E:=MBM^*$ is also a diagonal matrix; in particular,
$E=diag(b_1^2,\dots,b_n^2)$ where the $b_i^2$ are defined by
$a_i^2 + b_i^2 := 1$. Now we have (by the formula for the
quasi-projection operator):
\begin{align*}
M(\har{A}{B})M^* &= 2MAM^*(M(A+B)M^*)^{-1}MBM^* \\
   &= 2 diag(a_1^2 b_1^2,\dots,a_n^2 b_n^2) .
\end{align*}
Thus $M(\har{A}{B})M^*$ is positive semi-definite and hence $\har{A}{B}$ is
positive semi-definite.

\end{Proof}


\section{An iso-spectral flow which preserves zeros.}

As above, let $\Delta \subseteq \{1,2,\dots,n\} \times \{1,2,\dots,n\}$ be a
set
of pairs $(i,j)$ of indices which satisfy the following conditions
for all $i,j = 1,2,\dots,n$:
\begin{align*}
(i,i) &\in \Delta, \tag{nz1} \\
(i,j) \in \Delta &\text{ implies } (j,i) \in \Delta. \tag{nz2}
\end{align*}
Recall that we are using $Sym(n)$ to denote the vector space of
symmetric $n\times n$ matrices and we are using $Sym(\Delta)$ to
denote the subspace of $Sym(n)$ consisting of the symmetric matrices
which are zero outside of $\Delta$.
The set $\Delta$ of pairs of indices represents the nonzero pattern of
interest. The first condition on $\Delta$ implies that the diagonal
matrices are a subspace of $Sym(\Delta)$. The second condition simply
says that the pattern $\Delta$ is symmetric.  We want to consider some
iso-spectral flows on $Sym(\Delta)$.

We use $[X,Y]:= XY -YX $ to denote the {\bf commutator} of two
square matrices. Note that if $X$ is symmetric and $K$ is
skew-symmetric then $[X,K]$ is symmetric. Furthermore, we use
$O(n)$ to denote the orthogonal group. For a symmetric matrix $X$,
let
$$
\omega.X := O(n) \to Sym(n) : Q \mapsto QXQ^T .
$$
Then the image of $\omega.X$ is the iso-spectral surface,
$Iso(X)$, determined by $X$.
We can regard $\omega.X$ as a map from one manifold to
another. In particular we can differentiate this map at the identity $I$ to
obtain the  following linear map:
$$
D(\omega.X).I = Tan.O(n).I \to Tan.Sym(n).X: K \mapsto [K,X] .
$$
The space tangent to $O(n)$ at the identity $I$ may be
identified with the skew-symmetric matrices; in symbols,
$$
Tan.O(n).I = Skew(n) := \{K\in \R^{n \times n} : K^T = -K \}.
$$
(See, for example, Curtis~[1984].)
Clearly we can also identify $Tan.Sym(n).X$ with $Sym(n)$.
Hence, we define a map $l.X$ as a linear map from $Skew(n)$ to $Sym(n)$ by
$$
l.X := D(\omega.X).I= Skew(n) \to Sym(n) : K \mapsto [K,X] .
$$
It is not hard to prove that
the space tangent to $Iso(X)$ at $X$ is the image of the linear map
$D(\omega.X).I$; in symbols,
$$
Tan.Iso(X).X = \{ [K,X]: K\in Skew(n) \} .
$$
(For details see
Warner~[1983] chapter 3: Lie groups, section: homogeneous manifolds.)

\begin{Remark} Note that if $X$ has distinct eigenvalues then (by
the spectral theorem) the map $l.X$ is injective. However, if some
of the eigenvalues of $X$ are repeated then $l.X$ is not
injective. This is one of the reasons that we prefer to use
quasi-projection rather than projection. See the appendix which
compares projections and quasi-projections.
\end{Remark}

Recall that we are using the Frobenius inner product on $n$-by-$n$
matrices: $\langle X,Y\rangle:= Trace(X Y^T)$. We list a few properties of this
inner product: $\langle X Y,Z\rangle = \langle X,Z Y^T\rangle = \langle
Y,X^T Z\rangle$ and $\langle [X,Y],Z\rangle = 
\langle X,[Z,Y^T]\rangle = \langle Y,[X^T,Z]\rangle$.

The adjoint $(l.X)^*$ of $l.X$ is the following map:
$$
(l.X)^* = Sym(n) \to Skew(n): Y \mapsto [Y,X]
$$
since, for every symmetric matrix $Y$ and every skew-symmetric
matrix $K$, $\langle [K,X],Y \rangle = \langle K,[Y,X] \rangle$.
The composition of $l.X$ with its adjoint is a ``double bracket'':
$$
(l.X)(l.X)^* = Sym(n)\to Sym(n): Y \mapsto [[Y,X],X] .
$$
Note that for any $Y\in Sym(n)$, we have that $(l.X)(l.X)^*.Y$ is
tangent to the iso-spectral surface $Iso(X)$ at $X$.

We shall also use the map $m: Sym(n) \to Sym(\Delta)$ which is defined
as follows: For any symmetric matrix $Y$, let $m.Y$ denote the matrix
defined by $m.Y(i,j):= Y(i,j)$ if $(i,j)$ is in $\Delta$ and
$m.Y(i,j):=0$ if $(i,j)$ is not in $\Delta$. Note that $m$ is the
orthogonal projection of $Sym(n)$ onto $Sym(\Delta)$. In particular,
we have $m=m^*=m^2$.

We want to consider vector fields on $Sym(\Delta)$ which are
iso-spectral.  We can obtain such vector fields by quasi-projection.
Let $v:Sym(n)\to Sym(n)$ be any smooth map on $Sym(n)$. From $v$
we can obtain an iso-spectral vector field on $Sym(\Delta)$ by
quasi-projection as follows. For any symmetric matrix $X$, let $\rho.X
:=\har{A.X}{B.X}$ be the quasi-projection map determined by
$A.X:=(l.X)(l.X)^*$ and $B.X:=m$. Since these latter two linear maps are
positive semi-definite, the results of the last section apply here. We
shall use those results without explicitly citing particular
propositions.  In particular, note that for any symmetric matrix $Y$,
the symmetric matrix $\rho.X.Y$ is in the intersection of the range of
$(l.X)(l.X)^*$ and $m$; in symbols,
$$
\rho.X.Y \in Tan.Iso(X).X \cap Sym(\Delta) .
$$
We have the following iso-spectral vector field on
$Sym(\Delta)$:
$$
Sym(\Delta) \to Sym(\Delta): X \mapsto \rho.X(v.X) .
$$
The corresponding differential equation is $X^\prime = \rho.X(v.X)$.
We can rewrite this differential equation as a differential (linear)
algebraic equation as follows:
\begin{align*}
X^\prime &= (l.X)(l.X)^*(\lambda + v.X)\,, \\
(l.X)(l.X)^*( \lambda +v.X) &= m(-\lambda + v.X) .
\end{align*}
Note that the second of these equations is a linear equation for the
unknown symmetric matrix $\lambda$. The vector field is determined by
solving this second equation for $\lambda$ and substituting the
solution into the first equation.

Using the formulas for $l.X$ and $(l.X)^*$, we get
$$
(l.X)(l.X)^*(\lambda + v.X) = [[\lambda+v.X,X],X] .
$$
Substituting this simplification into the differential algebraic
equation, we get
\begin{align*}
X^\prime &= [[\lambda+v.X,X],X]\,, \\
[[\lambda + v.X,X],X] &= m(-\lambda + v.X) .
\end{align*}

We now turn our attention to a specific flow. This flow is
determined by the optimization problem (Problem~1) that we
mentioned in the introduction.
We shall see that we can solve this problem by finding a vector
field on $Sym(\Delta)$ associated with the objective function $f$
which is iso-spectral.
We obtain $X-D$ for the gradient of $f$ at $X$, in symbols $\nabla
f.X = X-D$.
We can get an iso-spectral vector field by orthogonal projection
of $\nabla f.X$ onto the intersection $Tan.Iso(X).X \cap
Sym(\Delta)$. We prefer to quasi-project instead. (We explain this
preference in an appendix.) We simply substitute the negative of
the gradient into the formulas given above. We get the following
system:
\begin{align*}
X^\prime &= [[\lambda + D,X],X], \tag{de1} \\
[[\lambda + D,X],X] &= -m(\lambda+X-D) . \tag{de2}
\end{align*}
We call the flow generated by this system the {\bf quasi-projected gradient
  flow} determined by the objective function $f$. We summarize the
properties of this flow in the following proposition.

\begin{Proposition}
Let $D$ be a symmetric matrix. Then the system (de1) and (de2) generating the
quasi-projected gradient flow  has the following properties:

\begin{itemize}
\item[(i)] The quasi-projected gradient flow preserves eigenvalues and the
       nonzero pattern of interest.

\item[(ii)] The function $f(X):=(1/2)\langle X-D,X-D \rangle$ is
  non-increasing along solutions of this system.

\item[(iii)] A point $E\in Sym(\Delta)$ is an equilibrium point of this system
  iff it satisfies the conditions 
  \begin{align*}
    [\lambda + D, E] &= 0, \text{ and } \tag{e1} \\
    m(\lambda + E -D) &= 0 . \tag{e2}
  \end{align*}
  for some symmetric matrix $\lambda$.

\item[(iv)] If a matrix $E\in Sym(\Delta)$ commutes with $D$ then $E$
  is an equilibrium point of this system.
\end{itemize}
\end{Proposition}
    
\begin{Proof}
(i) 
That this flow preserves eigenvalues and the nonzero
pattern of interest is clear from the discussion above. The vector
field was chosen to have these properties. In particular, the vector
field preserves the nonzero pattern because $X^\prime =
-m(\lambda+X-D)$ has the nonzero pattern of interest. Also the vector
field preserves eigenvalues because $X^\prime=[[\lambda+D,X],X]$ is
tangent to the iso-spectral surface $Iso(X)$ at $X$.

(ii) 
Let $X(t)$ be any solution of the differential equation. Then,
since the quasi-projection operator $\rho.X=\har{(l.X)(l.X)^*}{m}$ is
positive semi-definite, we have
\begin{align*}
(f(X))^\prime &= \langle \nabla f.X, X^\prime \rangle
= \langle \nabla f.X, \rho.X(-\nabla f.X) \rangle \\
&= - \langle \nabla f.X , \rho.X(\nabla f.X) \rangle \leq 0 .
\end{align*}

(iii) 
Let $E\in Sym(\Delta)$ satisfy conditions (e1) and (e2). Then clearly 
$[[\lambda + D,E],E] = [0,E] = 0$ and $E$ is an equilibrium point of the
system (de1, de2).
On the other hand, if $E\in Sym(\Delta)$ is an equilibrium point then
(de1) implies $[[\lambda+D,E],E]=0$ and (de2) implies (e2). We then also
get 
$$
0 = \langle [[\lambda + D,E],E] , \lambda + D \rangle
 = \langle [\lambda+D,E],[\lambda+D,E] \rangle ,
$$
which implies (e1).

(iv)
Take $\lambda:=D-E$. Then (e2) is trivially satisfied and for (e1) we
have
$$
[\lambda + D, E] = [2D-E,E] = 2[D,E] = 0 .
$$
\end{Proof}

\begin{Remark}
We should say a few words about convergence of this
system.  (We intend to discuss convergence more fully in a future
paper.)  Note that the map $\omega.A$ is a smooth map from $O(n)$ onto
$Iso(A)$. Hence $Iso(A)$ is compact since $O(n)$ is compact. From part
(i) of the proposition, we then see that every solution starting in the
iso-spectral surface $Iso(A)$ remains in that surface and is entire.
(In particular,``blowup'' is not possible.) Again using compactness, we
see that every such solution has $\omega$-limit points. If the
equilibrium points on the iso-spectral surface are isolated (which we
expect is usually true) then every solution that starts in the
iso-spectral surface tends to an equilibrium point.  
\end{Remark}

Note that if $D$ is a diagonal matrix with distinct diagonal entries
and $E$ is a diagonal matrix then $E$ commutes with $D$. It follows
from part (iv) of the theorem that $E$ is an equilibrium point of the
quasi-projected gradient flow determined by $D$. In 2001 we
conjectured that diagonal matrices were the only equilibrium points of
this flow. In 2005 Bryan Shader found a counterexample to that 
conjecture. Here is a counterexample.

\begin{Example}
%
Let $a$ and $b$ be real non-zero parameters and consider 
the non-diagonal, symmetric matrix
\[
E:=\begin{pmatrix}
       0 & a&  0\\
       a & 0&  b\\
       0 & b&  0
\end{pmatrix}
\]
The matrix $E$ has the distinct eigenvalues $0$ and $\pm\sqrt{a^2+b^2}$.

We set $\Delta$ as the non-zero pattern of $E$. 
We show, by suitably defining matrices $D$ and $\lambda$,
that $E$ is an equilibrium point of the quasi-projected 
gradient flow, i.e. satisfies conditions (e1) and (e2).

Let $y$ and $z$ be real parameters and take
$
\lambda +D := yE + zE^2\,.
$
This clearly gives $[\lambda + D, E] = 0$, i.e. condition (e1) is 
satisfied. 
Furthermore,
\[
m(\lambda+E-D) = m(\lambda+D+E -2D)
=m((y+1)E + z E^2)-2D
=(y+1)E+z \cdot diag(E^2)-2D\,,
\]
where
\[
E^2 =\begin{pmatrix}
       a^2 & 0      &  ab\\
       0   & a^2+b^2&  0\\
       ab  & 0      &  b^2
\end{pmatrix}
\quad\text{and}\quad
diag(E^2):=\begin{pmatrix}
       a^2 & 0      &  0\\
       0   & a^2+b^2&  0\\
       0   & 0      &  b^2
\end{pmatrix}\,.
\]
Now, by choosing $y=-1$ and defining $D:=\frac{1}{2}z\cdot diag(E^2)$, 
we arrive at $m(\lambda+E-D)=0$, i.e. condition (e2) is satisfied. 
Furthermore, if $z\ne 0$ and $|a|\ne|b|$ then $D$ has the required
distinct diagonal entries.

A numerical experiment
shows that the equilibrium point $E$ with $a:=1$ and $b:=2$ and $z:=2$
is not stable.
\end{Example}

We now conjecture that if $D$ is a diagonal matrix with distinct diagonal
entries then diagonal matrices are the only \emph{stable} 
equilibrium points of the 
quasi-projected gradient flow determined by $D$.

\begin{Remark}
A set $S$ in a topological space $T$ is called \emph{nowhere-dense} if 
the interior of its closure is empty. A set $S\subset T$ is called 
\emph{generic}  if it is open and dense.
Note that if $S$ is closed then it is nowhere-dense iff $T\backslash{}S$ 
is generic.

Let $V$ be a (finite-dimensional) vector space over the reals $R$. Let 
$f:V \to R$ be a real-valued function on $V$. Note that if $f(x)$ is a
polynomial in the components of $x$ with respect to some basis for $V$
then $f$ has this property for every choice of basis. In this case we
say that $f$ is a \emph{polynomial function}.

Proposition: Let $f:V \to R$ be a polynomial function. If $f$ is not the
zero polynomial then the variety, $Variety(f) := \{x \in V: f(x)=0\}$, of
$f$ is nowhere-dense.

Remark: We use the standard topology on $V$.

Proof: Note the variety is closed.  Suppose that the variety is not
nowhere-dense. Then $f$ vanishes on an open subset of $V$. It follows that
$f$ is identically $0$. \qed

%

Here is an application involving determinants. 

Example: Consider the determinant function $\det:R^{n \times n} \to R$.
The set $\{M \in R^{n \times n} : \det.M = 0\}$ is nowhere-dense and
closed. Hence the set of non-singular $n \times n$ matrices is generic. 

Let $V$ and $W$ be vector spaces and let $f:V \to W$ be a map. Then $f$ is a
\emph{polynomial map} if the components $f_i(x)$, for $i=1,...,dim.W$, with
respect to some basis for $W$ are polynomial functions. Note that the
composition of two polynomial maps is a polynomial map.

At the beginning of Section~4, we will introduce
the assumption that the map $(A.X + m) : Sym(n) \to Sym(n)$  
with $A.X=(l.X)(l.X)^*=[[\cdot,X],X]$ is invertible for given $X\in Sym(n)$. 
Here we show that this is generic 
behavior if $X$ has distinct eigenvalues.
%
%
Hence, consider the map $A.X := Sym(n) \to Lin(Sym(n) \to Sym(n)) : X \mapsto
(Y \mapsto [[Y,X],X])$.
%

Proposition. 
The set $\{X \in Sym(n) : Range(A.X+m) = Sym(n)\}$ is generic.

Proof: Consider the polynomial function $Sym(n) \to R$ defined by 
$X \mapsto \det(A.X + m)$. 
(Here $\det$ is regarded as a real-valued function on the space of
linear maps $Lin(Sym(n) \to Sym(n))$.  
Note that $A$ and $m$ are polynomial maps.)
We show below that this function is not the zero function.
Then we have that $\{X \in Sym(n) : \det(A.X+m) = 0\}$ is 
nowhere-dense. Since this set is also closed we have that
\[
\{X \!\in Sym(n) : Range(A.X+m) = Sym(n)\}
= Sym(n) \backslash \{X \!\in Sym(n) : \det(A.X+m)\! = 0\}\,.
\]
is generic.
To complete the proof, we show that the function $X \mapsto \det(A.X + m)$ 
is not the zero function. Take $X=D:=diag(d_1,... d_n)$ where the $d_i$ 
are distinct. Then $A.D.Y(i,j) = (d_i - d_j)^2 Y(i,j)$. Furthermore, the
range of $m$ includes the diagonal matrices. These properties together show
that $Range(A.D + m) = Sym(n)$ and hence $\det(A.D + m)\ne 0$.\quad\qed

\noindent\mbox{}%
\end{Remark}

\section{Numerical results.}

We have implemented the quasi-projected gradient flow in a
Matlab program. This flow is iso-spectral and preserves zeros as discussed
in the previous section. 
In our implementation we assume that $Range((l.X)(l.X)^*+m)=Sym(n)$. 
We solve numerically for $t>0$ the initial value problem for
$X(t)$ given by
\[
X'=g(X):=2 m\left( (l.X)(l.X)^*+m \right)^{-1}
(l.X)(l.X)^*. (D-X),\quad X(0)=X_0\,,
\]
where $X_0\in Sym(\Delta)$ ($\Delta$ is defined by the nonzero
pattern of $X_0$ and kept constant) and $D$ is the diagonal matrix,
$D:=diag(1,2,\ldots,n)$. We refer to this flow as the {\bf Zero flow} in the
discussion of the examples and in the figures below.

The assumption on the ranges of $(l.X)(l.X)^*$ and $m$ guarantees the existence
of the inverse in the right-hand side of the differential
equation. This assumption is not satisfied in general as the following
example demonstrates.

\begin{Example}
Let $X$ be the circulant matrix with $-2$ on the diagonal and $1$ on
the first sub- and super-diagonal (and the corresponding corner
entries). The pattern $\Delta$ is defined as the nonzero pattern of
$X$. Now let $Y$ be any circulant matrix with nonzero pattern
completely outside of $\Delta$, i.e. $m.Y=0$. If $n=4$
this is, for instance, achieved by selecting $Y$ as the matrix with ones
on second sub- and super-diagonal. Since circulant matrices commute with each
other, and by the choice of the nonzero pattern of $Y$ we have
$((l.X)(l.X)^*+m)Y=0$. Thus $Y$ is a non-trivial element in the kernel
of the map and hence the inverse does not exist.
\end{Example}

By construction of the flow, the matrix $g(X(t))\in Sym(\Delta)$ for
all $t\ge 0$ and
when integrating the differential equation we ignore all matrix
elements outside the pattern $\Delta$ (these remain zero for all
$t>0$). Therefore the dimension of our differential equation is
reduced to the cardinality of $\Delta$ which may be significant less
than $n^2$. (We have currently not taken into account the symmetry of
the matrices.) However, we remark that we obtain intermediate matrices,
when evaluating the expression for $g(X(t))$ from the right to the
left, which can have nonzero entries outside of $\Delta$.

For the numerical solution of the initial value problem we employ
Matlab's explicit Runge-Kutta method of order 4(5) ({\tt rk45}) with
absolute and relative tolerance requirement set to $10^{-13}$. These
very stringent accuracy requirements reflect the fact that we are
currently interested in very accurate solutions to the initial value
problem and not (yet) in competitive numerical schemes for the
solution of sparse eigenvalue problems. Therefore, the cost of the numerical
computations are not considered in the following.

During the course of integration we monitor two characteristic
quantities of the flow.
\begin{enumerate}
\item The relative departure of the matrix $X(t)$ from the
iso-spectral surface associated with the initial matrix $X_0$. In
particular, we define
\[
d_{ev}(t):=\frac{\|ev(X_0)-ev(X(t))\|}{\|ev(X_0)\|}\,,
\]
where $ev(X)$ is the vector of sorted eigenvalues of $X$.
This quantity measures the quality of the time integration and
should be approximately constant in time and near the
machine accuracy ($10^{-14}$).
\item The relative size of the off-diagonal elements of $X(t)$ (with
respect to $X_0$). In particular, we define
\[
d_{off}(t):=\frac{\|X(t)-diag(X(t))\|_F}
{\|X_0-diag(X_0)\|_F}\,,
\]
where $diag(X)$ is the matrix containing the diagonal part of $X$ and
$\|\cdot\|_F$ is the Frobenius norm. This quantity measures the
convergence of the flow to a diagonal steady state and in conjunction
with a constant value of $d_{ev}(t)$ the convergence to the diagonal
matrix with elements corresponding to the eigenvalues of $X_0$.
\end{enumerate}

We compare the Zero flow with the ``double-bracket (DB) flow''. That is,
we also numerically solve the initial value problem
\[
X'=h(X):= [ [D, X], X],\quad X(0)=X_0\,,
\]
where $X_0\in Sym(n)$ and $D$ is the same diagonal matrix as
above. This flow is also iso-spectral and converges to a diagonal
matrix steady state with the eigenvalues of $X_0$ on the diagonal.
(See the appendix on the Toda flow and/or Driessel~[2004].) It
does not preserve the zero pattern of the initial matrix $X_0$ and
considerable fill-in can appear.
The double-bracket flow coincides with the Toda flow
if $X_0$ is a tridiagonal matrix.

We consider three different kinds of initial data in the next three
subsections.

\subsection{Example 1}

Here the initial value $X_0$ is a symmetric, tridiagonal
random matrix of dimension $6$:
{\small
\[
X_0:=\left(
\begin{array}{llllll}
    0.87& 1.23& 0   & 0   & 0   & 0\\
    1.23& 1.67& 0.62& 0   & 0   & 0\\
    0   & 0.62& 0.25& 1.17& 0   & 0\\
    0   & 0   & 1.17& 0.79& 1.87& 0\\
    0   & 0   & 0   & 1.87& 1.92& 1.63\\
    0   & 0   & 0   & 0   & 1.63& 1.8
\end{array}
\right)\,.
\]
}
We note that the DB flow preserves the tridiagonal pattern but
we do not exploit this fact in our implementation.

\begin{figure}[h]
\mbox{}\hfill
\begin{minipage}{0.6\hsize}
\includegraphics[width=\hsize]{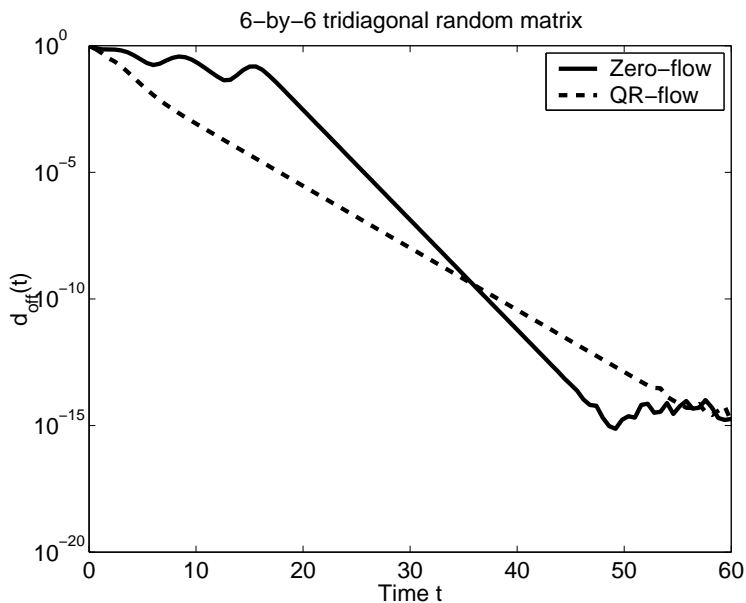}
\end{minipage}
\hfill\mbox{} \caption{\label{FigExample1Convergence} Convergence
history of the off-diagonal elements of the solution of Example~1
for the Zero and the DB flow.}
\end{figure}

We simulate the solution with this initial value until $t=60$ for
both flows. The maximum value of $d_{ev}(t)\approx 7\cdot
10^{-14}$ for the Zero flow and $\approx 2\cdot 10^{-14}$ for the
DB flow. This shows that for both flows the eigenvalues of the
initial matrix are preserved up to machine accuracy in the
numerical solution. In Figure~\ref{FigExample1Convergence}, we
plot the monitored values of $d_{off}(t)$ for both flows. We
observe that both converge to zero and that this happens slightly
faster for the DB flow initially but later the Zero flow converges
faster and reaches machine accuracy before the DB flow. The
results of this example show that for tridiagonal matrices the
Toda flow is different than our zero flow.

\subsection{Example 2}
\begin{figure}[h]
\mbox{}\hfill
\begin{minipage}{0.6\hsize}
\includegraphics[width=\hsize]{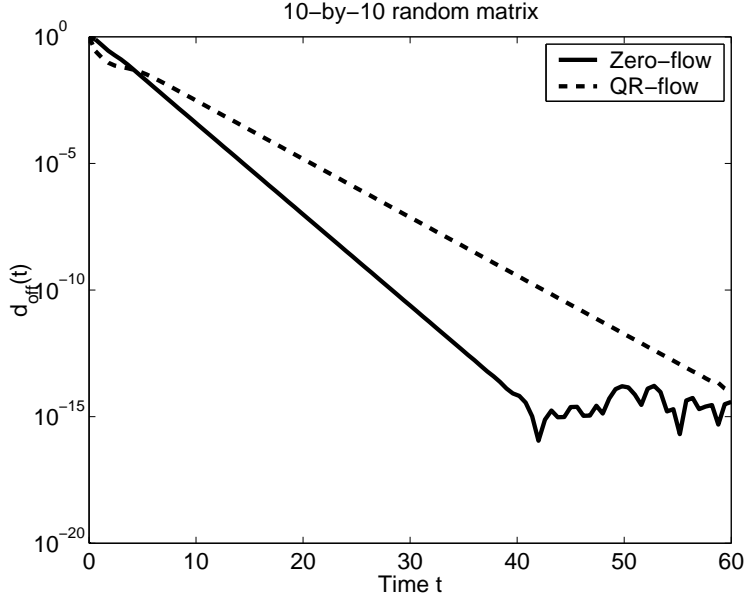}
\end{minipage}
\hfill\mbox{} \caption{\label{FigExample2Convergence} Convergence
history of the off-diagonal elements of the solution of Example~2
for the Zero and the DB flow.}
\end{figure}

In this example the initial value $X_0$ is a symmetric random matrix
of dimension $10$ with a random zero pattern:
{\small
\[
X_0:=\left(
\begin{array}{llllllllll}
    1.7&  0   &    0   & 0   & 0   & 0   & 1.92& 0   & 0.48& 1.25\\
    0  &  1.16&    1.16& 0.91& 1.56& 0   & 0   & 1.69& 0   & 0\\
    0  &  1.16&    0.48& 0   & 0.90& 0   & 0   & 0   & 0   & 0\\
    0  &  0.91&    0   & 0.66& 0.88& 0   & 0.93& 1.25& 0   & 1.39\\
    0  &  1.56&    0.9 & 0.88& 0.3 & 0   & 0   & 0   & 0   & 0\\
    0  &  0   &    0   & 0   & 0   & 0.94& 1.49& 0.37& 0.88& 0\\
    1.92& 0   &    0   & 0.93& 0   & 1.49& 1.12& 0.67& 0.4 & 0\\
    0   & 1.69&    0   & 1.25& 0   & 0.37& 0.67& 1.1 & 0   & 1.54\\
    0.48& 0   &    0   & 0   & 0   & 0.88& 0.4 & 0   & 0.44& 1.05\\
    1.25& 0   &    0   & 1.39& 0   & 0   & 0   & 1.54& 1.05& 1.2
\end{array}
\right)\,.
\]
} We simulate the solution with this initial value until $t=60$
for both flows. The maximum value of $d_{ev}(t)\approx 2\cdot
10^{-14}$ for the Zero flow and $\approx 6\cdot 10^{-15}$ for the
DB flow. This shows that for both flows the eigenvalues of the
initial matrix are preserved up to machine accuracy in the
numerical solution. In Figure~\ref{FigExample2Convergence}, we
plot the monitored values of $d_{off}(t)$ for both flows. We
observe that both converge to zero and that this happens slightly
faster for the Zero flow.

\subsection{Example 3}

\begin{figure}
\begin{minipage}{0.48\hsize}
\includegraphics[width=\hsize]{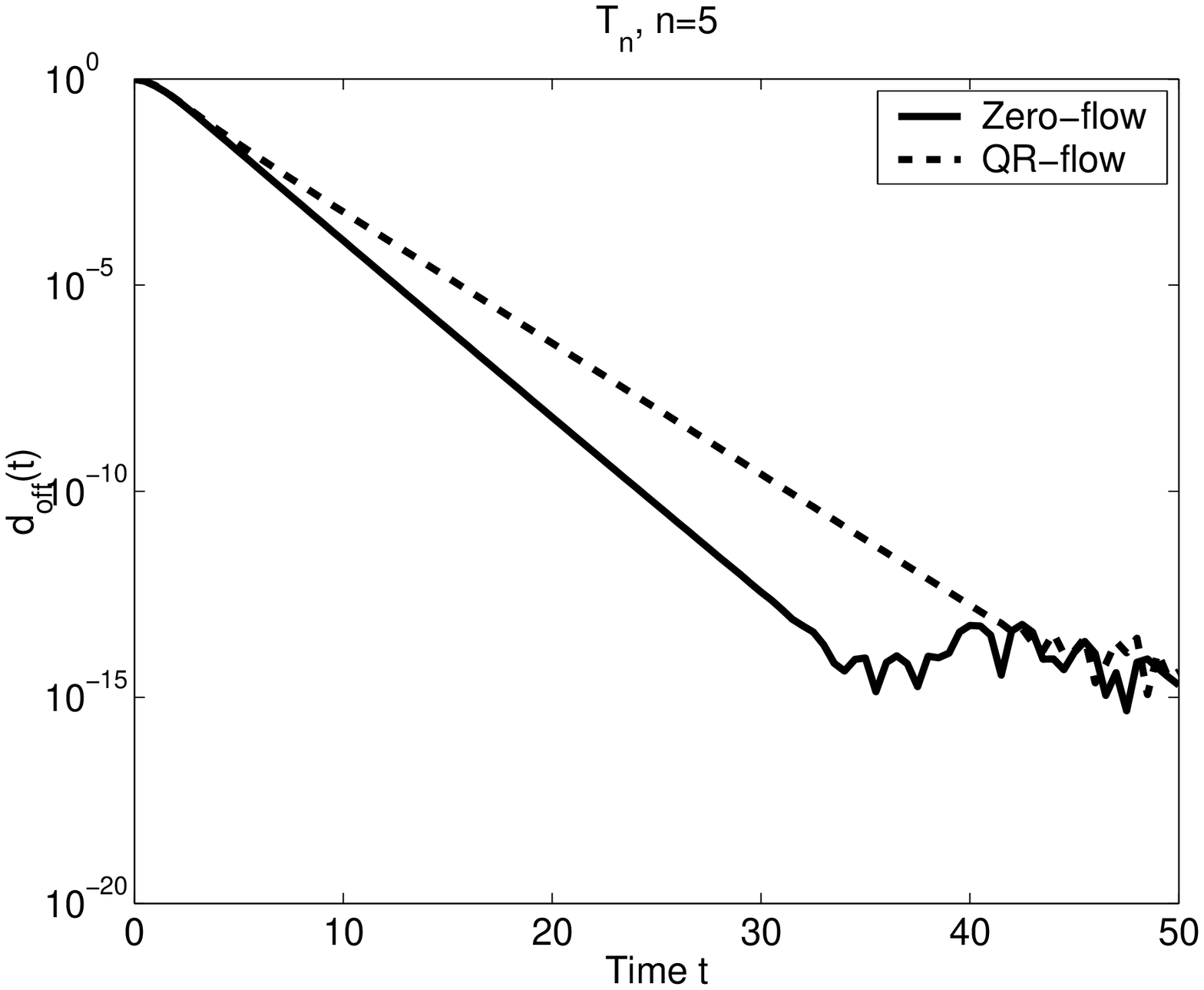}
\end{minipage}
\begin{minipage}{0.48\hsize}
\includegraphics[width=\hsize]{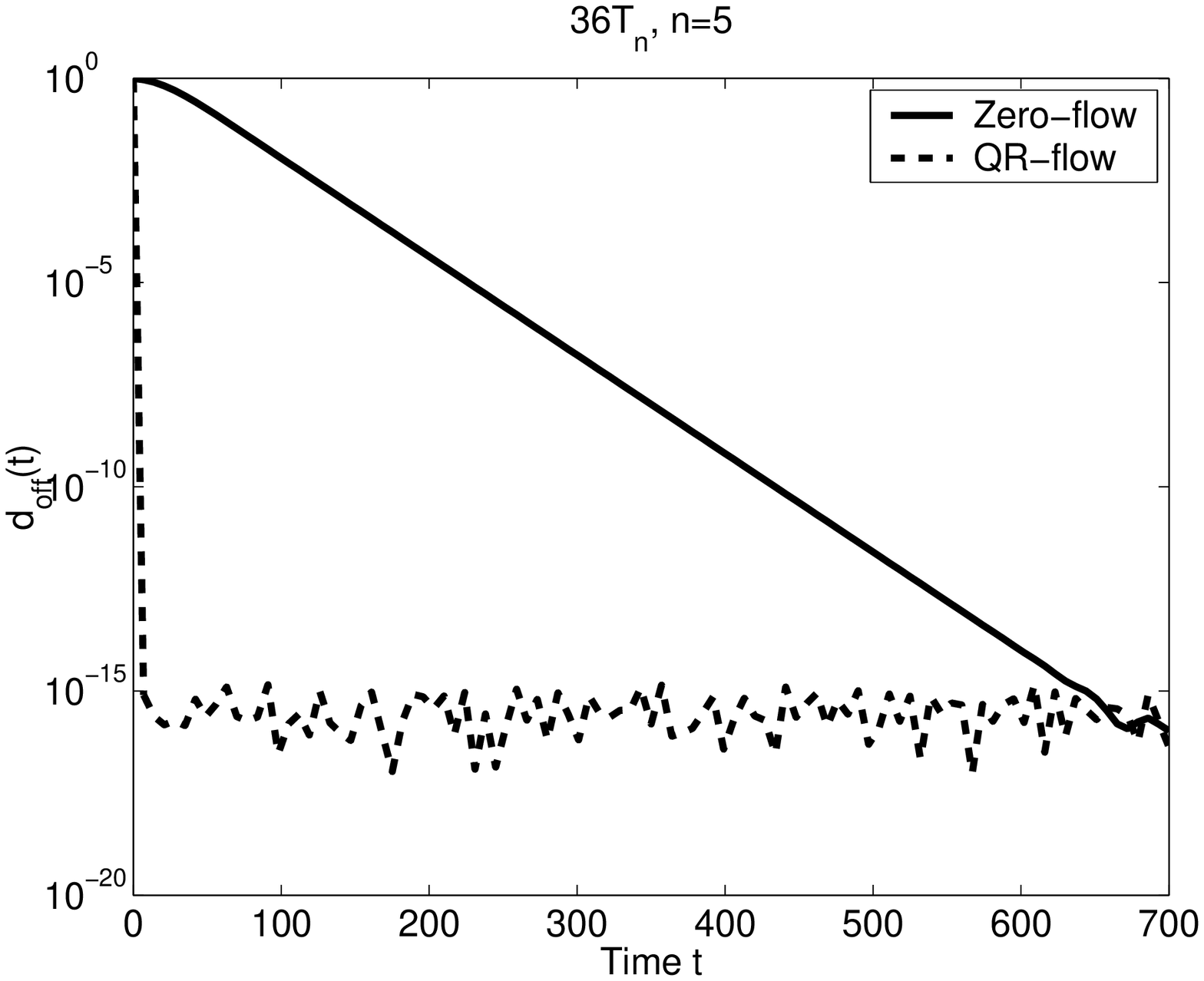}
\end{minipage}

\begin{minipage}{0.48\hsize}
\includegraphics[width=\hsize]{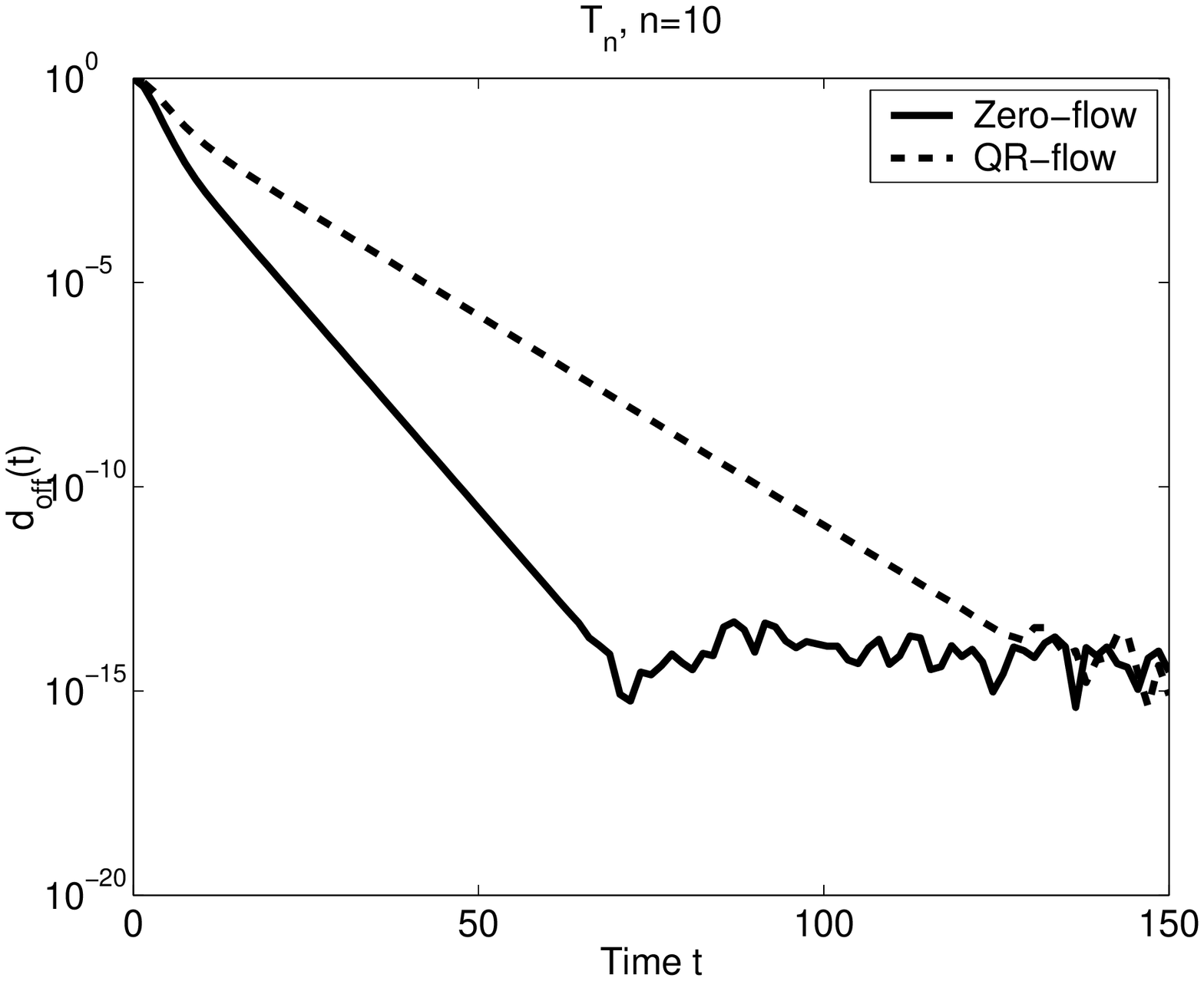}
\end{minipage}
\begin{minipage}{0.48\hsize}
\includegraphics[width=\hsize]{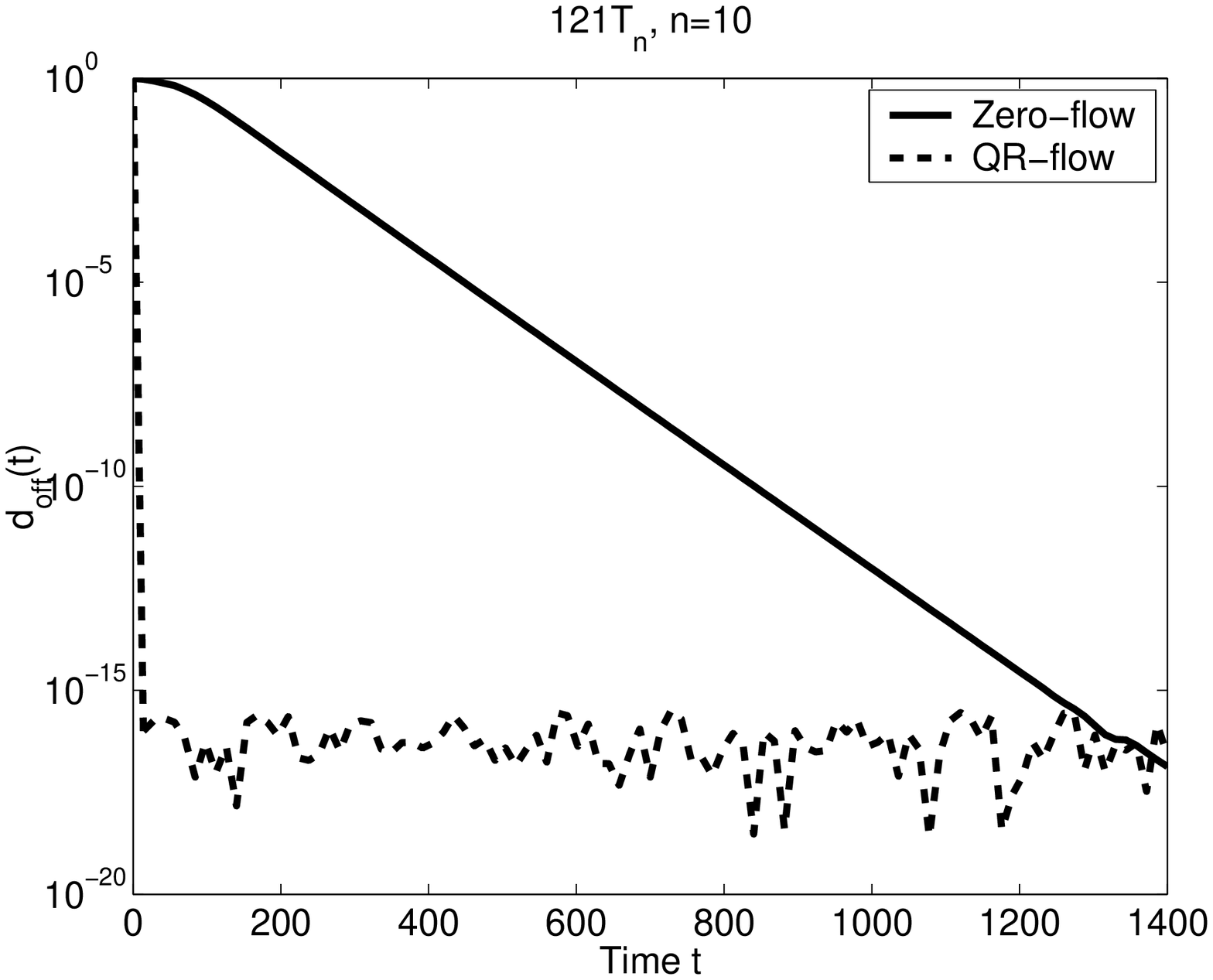}
\end{minipage}

\caption{\label{FigExample3Convergence} Convergence history of the
off-diagonal elements of the solutions of Example~3 for the Zero
and the DB flow for initial conditions $X_0=T_5$ (top left),
$X_0=\tilde T_5 = 36 T_5$ (top right), $X_0=T_{10}$ (bottom left),
and $X_0=\tilde T_{10} = 121 T_{10}$ (top right). }
\end{figure}

In the third example we consider tridiagonal matrices which arise when
one discretizes the boundary value problem  $u_{xx}=0, u(0)=u(1)=0$ by
standard second-order central differences.
Let $T_n:=tridiag(1,-2,1)\in Sym(n)$ and $\tilde T_n:=(n+1)^2 T_n$. Hence
$\tilde T_n$ corresponds to the discretization matrix of the boundary value
problem on an equidistant grid  with grid width $h:=1/(n+1)$.
The eigenvalues of both, $T_n$ and $\tilde T_n$, are distinct and
negative. We present results for the four cases $X_0=T_5, T_{10},
\tilde T_5$, and $\tilde T_{10}$ in Figure~\ref{FigExample3Convergence}.

We run these experiments to different final times as can be seen
in the plots. We note that the values of $d_{ev}(t)$ are in the
range $10^{-15}$ to $10^{-13}$ for all values of $t$ considered.
Again this demonstrates that the numerical solution does only
insignificantly drift off the iso-spectral surface associated with
the initial matrix. Both the Zero flow and the DB flow converge to
the diagonal matrix containing the eigenvalues of the initial
data. However, whereas the Zero flow does so much faster than the
DB flow for the matrices $T_n$, the situation is the opposite for
the scaled matrices $\tilde T_n$. The change in the convergence
speed of the DB flow for different initial matrices $T_n$ and
$\tilde T_n$ is precisely explained by the following proposition.

\begin{Proposition}
If $X(t)$ is the solution of the double-bracket flow with initial value
$X_0$ then $c X(ct)$ is the solution of the double-bracket flow with
initial value $c X_0$, $c>0$.
\end{Proposition}
\begin{Proof}
We can write $h(cX)=[[D, cX], cX] = c^2 [[D, X], X] = c^2 h(X)$. This
relation gives the desired scaling result.
\end{Proof}

We have not analyzed how scaling affects the Zero flow.

\appendix
\section{Comparing Projections with Quasi-Projections.}

Recall the following well-known result. (See, for example, Leon~[1986] Section
5.5:
Least-squares problems or Strang~[1980] Section 3.2: Projections onto subspaces
and least-squares approximation.)

\begin{Proposition} \textbf{(Least squares approximation).}
Let $L:U\to W$ be a linear map from one inner product space to
another. If $L$ is injective then, for any $b\in W$, the ``normal
equation'' $L^*Lx=L^*b$ has a unique solution which is given by
$(L^*L)^{-1}L^*b$. Furthermore, the linear map
$P:=L(L^*L)^{-1}L^*$ on $W$ has the following properties:
\begin{enumerate}
\item The range of $P$ equals the range of $L$: $Range.P=Range.L$.
\item The kernel of $P$ equals the orthogonal complement of the range of $L$:
$Kernel.P=(Range.L)^\bot$.
\item The map $P$ is the projection of $W$ on $Range.L$ along
$(Range.L)^\bot$ which corresponds to the decomposition
$W=Range.L\oplus (Range.L)^\bot$. In particular, $P^2=P$ and $P^*=P$.
\end{enumerate}
\end{Proposition}

The map $L^+=(L^*L)^{-1}L^*:W\to U$ is the Moore-Penrose pseudo-inverse
of $L$. The map $P=LL^+$ is the \textbf{projection map associated with
the least squares problem} $Lx\sim b$. The projected vector $Pb$ is the
element of the range of $L$ which is closest to $b$ in the least
squares sense.

Driessel~[2004] observed the following: It is often difficult to
directly use the projection map $P$. If $L$ is not injective then the
inverse of $L^*L$ does not exist. Even when $L$ is injective it is
often difficult to compute $(L^*L)^{-1}$ - for example, if $L$ is
ill-conditioned or the dimension of the vector space is large. We can
often avoid these difficulties by using the linear map $LL^*:W\to W$
instead of the projection map $P$. For the map $L^*L$ we have the following
analogue of the last proposition.

\begin{Proposition} \textbf{(Quasi-projection)}
Let $L:U\to W$ be a linear map between two inner product spaces. Then
the map $LL^*:W\to W$ has the following properties:
\begin{enumerate}
\item The range of $LL^*$ equals the range of $L$: $Range(LL^*)=Range.L$.
\item The kernel of $LL^*$ equals the orthogonal complement of the
range of $L$: $Kernel(LL^*)=(Range.L)^\bot$.
\item The map $LL^*$ is positive semi-definite.
\end{enumerate}
\end{Proposition}
We include the proof of this proposition from Driessel~[2004] for
completeness.

\begin{Proof} Here is the proof of the first assertion. It is obvious that
$Range(LL^*)\subseteq Range.L$. We want to see the other
inclusion. Consider any element $Lx$ in the range of $L$. Let $x=y+z$
where $y\in (Kernel.L)^\bot$ and $z\in Kernel.L$. Since
$Range.L^*=(Kernel.L)^\bot$ we have $Lx=Ly$ is an element of
$Range(LL^*)$. Here is the proof of the third assertion:
$$
\langle u,LL^*v \rangle = \langle L^*u,L^*v \rangle =
\langle L^{**}L^*u,v \rangle = \langle LL^*u,v \rangle
$$
since $L^{**}=L$. Finally we consider the second assertion. By the first
and third assertions we have
$$
(Range.L)^\bot=(Range(LL^*))^\bot=Kernel(LL^*)^*=Kernel(LL^*).
$$
\end{Proof}

Driessel~[2004] called the map $LL^*:W\to W$ the
\textbf{quasi-projection map associated with the least squares
problem} $Lx\sim b$.  Driessel~[2004] also compared the projection $P$
with the quasi-projection $LL^*$ as follows: Since the restriction
$LL^*: Range.L\to Range.L$ of $LL^*$ is self-adjoint, we can find a
basis of $Range.L$ consisting of eigenvectors: $LL^*w_i=\lambda_i w_i$
for $i=1,2,\dots,m$ where $m$ is the dimension of $Range.L$. For any
$w \in W$ let $w=r+s$ where $r\in Range.L$ and $s\in
(Range.L)^\bot$. We have $Pw=r=\Sigma\langle r,w_i \rangle w_i$ and
$LL^*w=LL^*r=\Sigma\langle r,w_i \rangle \lambda_i w_i$. Thus $LL^*$
is a projection followed by an eigenvalue-eigenvector scaling.  (Also
note that $LL^*P=PLL^*=LL^*$.) Note that each $\lambda_i$ is
non-negative. It follows that the signature of $P$ is the same as the
signature of $LL^*$. We regard congruence as the appropriate geometry
for the study of quasi-projections. In summary, we regard the use of
the quasi-projection operators $LL^*$ as simpler, more direct and more
robust than the use of the projection operator $P$.

We want to establish propositions like the last two for a pair of
linear maps.  Let $U$, $V$ and $W$ be finite-dimensional inner product
spaces and let $L:U\to W$ and $M:V\to W$ be linear maps. We consider
the following problem:
\begin{Problem}\textbf{(Projection)} Given a vector $c\in W$, find the
vector $\hat c \in W$ which is in the intersection $Range.L\cap Range.M$
and is closest to $c$.
\end{Problem}
We can formulate this problem as a constrained optimization problem as
follows. Let
\begin{eqnarray*}
f.c &:=& U\times V \to \R :
(x,y) \mapsto (1/2)(\langle Lx-c,Lx-c \rangle + \langle My-c,My-c \rangle) , \\
k &:=& U\times V \to W: (x,y) \mapsto Lx-My .
\end{eqnarray*}
\begin{Problem} \textbf{(Optimization)} Given $c\in W$, find the pair in
$U\times V$ which minimizes $f.c(x,y)$ subject to the constraint
$k(x,y)=0$.
\end{Problem}
If $(u,v)$ is the solution of this optimization problem then
$\hat c=Lu =Mv$ is the solution of the projection problem.

We begin our analysis of the optimization problem by computing the
derivative of the objective function $f.c$. We have
$$
D(f.c)(x,y)(dx,dy)= \langle Lx-c,L\, dx\rangle + \langle My-c,M\, dy\rangle
= \langle L^*(Lx-c),dx\rangle + \langle M^*(My-c),dy \rangle .
$$
We use the standard Cartesian inner product on $U\times V$; that is,
for $(x,y)$ and $(x^\prime,y^\prime)$ in $U\times V$, we take
$\langle (x,y),(x^\prime,y^\prime) \rangle :=
\langle x,x^\prime \rangle + \langle y,y^\prime \rangle $.
From the equation for the derivative of $f.c$, we easily recognize the
gradient of $f.c$:
$$
\nabla(f.c)(x,y) = (L^*(Lx-c),M^*(My-c)) .
$$

Let $(u,v)\in U\times V$ be the solution of the optimization problem.
By the well-known Lagrange multiplier theorem, we have the condition
$$
\nabla(f.c)(u,v) \in (Kernel(Dk(u,v)))^\bot = Range(Dk(u,v))^* .
$$
Since $k$ is linear, we have $Dk(u,v)=k$. Next we
compute $k^*$: For $x\in U$, $y\in V$, and $z\in W$ we have
\begin{eqnarray*}
\langle k(x,y),z \rangle &=& \langle Lx-My,z \rangle
= \langle Lx,z \rangle - \langle My,z \rangle \\
&=& \langle x,L^*z \rangle - \langle y,M^*z \rangle =
\langle (x,y),(L^*z,-M^*z) \rangle .
\end{eqnarray*}
In other words,
$
k^* = W \to U\times V : z\mapsto (L^* z,-M^*z) .
$
Hence the Lagrange condition $\nabla(f.c)(u,v)\in Range.k^*$ is equivalent
to the following one:
$$
\exists \lambda \in W, L^*(Lu-c)=L^*\lambda \text{ and } M^*(Mv-c) =
-M^*\lambda .
$$
Adding the constraint condition, we get the following system of (linear)
equations for $\lambda$,$u$ and $v$:
\begin{eqnarray*}
L^*Lu - L^*\lambda &=& L^*c, \\
M^*Mv + M^*\lambda &=& M^*c, \\
Lu &=& Mv .
\end{eqnarray*}
We now assume that $L$ and $M$ are injective. Then $L^*L$ and $M^*M$ are
invertible. We can apply (block) Gaussian elimination to the last system
of equations; we get:
\begin{eqnarray*}
u - L^+\lambda &=& L^+c, \\
v + M^+\lambda &=& M^+c, \\
(LL^+ +  MM^+)\lambda &=& (-LL^+ +  MM^+)c
\end{eqnarray*}
where $L^+=(L^*L)^{-1}L^*$ and $M^+=(M^*M)^{-1}M^*$ are the Moore-Penrose
pseudo-inverses of $L$ and $M$ respectively. This last set of equations
implies the following set:
\begin{eqnarray*}
Lu - P\lambda &=& Pc, \\
Mv + Q\lambda &=& Qc, \\
(P + Q)\lambda &=& (-P + Q)c
\end{eqnarray*}
where $P:=LL^+$ and $Q:=MM^+$. Note that $P$ and $Q$ are the
orthogonal projections of $W$ onto $Range.L$ and $Range.M$ respectively.

For any vector $c$ in $W$, we are led to consider the following system
of linear equations for $w$ and $\lambda$ in $W$:
\begin{align*}
w - P\lambda &= Pc, \tag{p1}\\
(P + Q)\lambda &= (-P + Q)c \tag{p2}
\end{align*}
(We get this set of equations from the preceding set by setting
$w:=Lu=Mv$
and then omitting the redundant second equation.) Note that these are the
quasi-projection equations determined by $P$,$Q$, and $c$. Since $P$ and
$Q$ are positive semi-definite, the results in the section on
quasi-projections apply. In particular, we have the following corollaries.
\begin{Corollary}\textbf{(Uniqueness)}
For any $c$ in $W$, if $(w_1,\lambda_1)$ and $(w_2,\lambda_2)$ are
solutions of the equations (p1) and (p2) then $w_1=w_2$,
$P\lambda_1=P\lambda_2$ and $Q\lambda_1=Q\lambda_2$.
\end{Corollary}
Note that $\lambda$ is uniquely determined iff $P+Q$ is surjective.
\begin{Corollary}\textbf{(Existence)}
For all $c$ in $W$ there exist $w$ and $\lambda$ in $W$ satisfying
(p1) and (p2).
\end{Corollary}

As in the section on quasi-projections, we use $\har{P}{Q}$ to denote the
linear operator on $W$ which maps a vector $c$ to the unique vector $w$
which satisfies the following condition: There exists $\lambda$ in $W$
such that the pair $(w,\lambda)$ satisfies equations (p1) and (p2).
\begin{Corollary}\textbf{(Quasi-Projection Formulas)} The quasi-projection
operator $\har{P}{Q}$ satisfies
$$
\har{P}{Q}=2P(P+Q)^+Q=2Q(P+Q)^+P .
$$
Furthermore, $\har{P}{Q}$ is the ortho-projection of $W$ on $Range.P\cap Range.Q$.
\end{Corollary}
\begin{Proof}
Claim: If $c\in Range.P\cap Range.Q$ then $\har{P}{Q}c=c$.

We have $Pc=Qc=c$. It follows that taking $w:=c$ and $\lambda:=0$ gives us
a solution of (p1) and (p2).
\end{Proof}

In summary we have the following analogue of the proposition concerning
least squares approximation involving a single linear map.

\begin{Proposition} Let $U$, $V$ and $W$ be inner product spaces and let
$L:U\to W$ and $M:V\to W$ be injective linear maps. Let $P:=LL^+$ and
$Q:= MM^+$. Then the map $\har{P}{Q}$ has the following properties:
\begin{enumerate}
\item The range of $\har{P}{Q}$ equals the intersection of the ranges
of $L$ and $M$: $Range.\har{P}{Q} = Range.L \cap Range.M$.
\item The kernel of $\har{P}{Q}$ equals the orthogonal complement of the
intersection of the ranges of $L$ and $M$:
$Kernel.\har{P}{Q} = (Range.L\cap Range.M)^\bot$.
\item The map $\har{P}{Q}$ is the projection of $W$ onto $Range.L\cap Range.M$
along $(Range.L\cap Range.M)^\bot$ which corresponds to the decomposition
$$
W=(Range.L\cap Range.M)\oplus (Range.L\cap Range.M)^\bot .
$$
\end{enumerate}
\end{Proposition}

We also want to establish an analogue of the proposition concerning
the quasi-projection associated with a single linear map. We do so by
setting $A:=LL^*$ and $B:=MM^*$. We then consider the quasi-projection
equations determined by $c\in W$ and these maps:
\begin{align*}
w - A \lambda &= Ac , \tag{q1} \\
(A+B) \lambda &= (B-A)c . &\tag{q2}
\end{align*}
Note that $A$ and $B$ are positive semi-definite. Hence the results in
the section on quasi-projections apply. In particular, we have the
following result.

\begin{Proposition} Let $U$, $V$ and $W$ be inner product spaces and let
$L:U\to W$ and $M:V\to W$ be linear maps. Let $A:=LL^*$ and
$B:= MM^*$. Then the map $\har{A}{B}$ has the following properties:
\begin{enumerate}
\item The range of $\har{A}{B}$ equals the intersection of the ranges
of $L$ and $M$: $Range.\har{A}{B} = Range.L \cap Range.M$.
\item The kernel of $\har{A}{B}$ equals the orthogonal complement of the
intersection of the ranges of $L$ and $M$:
$Kernel.\har{A}{B} = (Range.L\cap Range.M)^\bot$.
\item The map $\har{A}{B}$ is positive semi-definite.
\end{enumerate}
\end{Proposition}

We regard the use of the quasi-projection operator $\har{A}{B}$ as
simpler, more direct and more robust than the use of the projection
operator $\har{P}{Q}$.  In particular, we do not need to compute
$(L^*L)^{-1}$ and $(M^*M)^{-1}$ when using $\har{A}{B}$. The signature
of $\har{A}{B}$ was determined in the proof of the proposition
concerning quasi-projection formulas in the section on
quasi-projections. It is easy to see that $\har{A}{B}$ and $\har{P}{Q}$
have the same signature.

\section{On the Toda flow.}

The flows that we describe are related to the QR algorithm and the
Toda flow. For a square matrix $X$ let $X_l, X_d$ and $X_u$ denote the
strictly lower triangular, diagonal and strictly upper triangular part
of $X$. The {\bf Toda flow} or {\bf QR flow} is the flow associated
with the following differential equation:
$$
X^\prime = [X, X_l - X_l^T] .
$$
We use $[X,Y]:= XY -YX $ to denote the {\bf commutator} of two
square matrices. Note that $X_l - X_l^T$ is skew-symmetric. Also
note that if $X$ is symmetric and $K$ is skew-symmetric then
$[X,K]$ is symmetric. Hence we can (and shall) view the Toda flow
as a dynamical system in the space of symmetric matrices.

For symmetric matrix $X$, consider the following map determined by
$X$:
$$
\omega.X := O(n) \to Sym(n): Q \mapsto QXQ^T .
$$
Note that the image of this map is the iso-spectral surface $Iso(X)$. We
differentiate $\omega.X$ to get the following linear map:
$$
D(\omega.X).I = Tan.O(n).I \to Tan.Sym(n).X: K \mapsto [K,X] .
$$
Recall that the space tangent to $O(n)$ at the identity $I$ may be
identified with the skew-symmetric matrices; in symbols,
$$
Tan.O(n).I = Skew(n) := \{K\in \R^{n \times n} : K^T = -K \}.
$$
(See, for example, Curtis~[1984].)  Clearly we can also identify
$Tan.Sym(n).X$ with $Sym(n)$. We shall regard $D(\omega.X).I$ as a
linear map from $Skew(n)$ to $Sym(n)$. It is not hard to prove that
the space tangent to $Iso(X)$ at $X$ is the image of the linear map
$D(\omega.X).I$; in symbols,
$$
Tan.Iso(X).X = \{ [K,X]: K\in Skew(n) \} .
$$
(For details see Warner~[1983]
chapter 3: Lie groups, section: homogeneous manifolds.)
Since the vector field $X \mapsto [X,X_l - X_l^T]$ of the Toda flow is
tangent to $Iso(X)$, it follows that the Toda flow is iso-spectral,
that is, it preserves eigenvalues.  It is well-known that the Toda
flow is iso-spectral; for details, see, for example, Demmel~[1997]
Section 5.5: ``Differential Equations and Eigenvalue Problems'' and
the references there. The relationship between the Toda flow and the
QR algorithm is also fairly well-known; again see, for example,
Demmel~[1997].

The QR algorithm and the Toda flow do have limited zero-preserving
properties. We say that a symmetric pattern of interest $\Delta$ is a
{\bf staircase pattern} if $\Delta$ is ``filled in toward the
diagonal'', that is, for all $i<j$, if $(i,j)$ is in $\Delta$ then so
are $(i,j-1)$ and $(i+1,j)$.
Arbenz and Golub~[1995] showed that the QR algorithm
preserves symmetric staircase patterns and only such sparseness.
Ashlock, Driessel and Hentzel~[1997a] showed that the Toda flow
preserves symmetric staircase patterns and only such sparseness.  Here
we aim to preserve arbitrary sparseness.

\begin{Remark}
For an earlier attempt to
generalize the Toda flow to other zero-preserving flows, see Ashlock,
Driessel and Hentzel~[1997b]. This attempt had only very limited
success. Chu and Norris~[1988] designed flows on the symmetric
matrices which converge to $Sym(\Delta)$ for various $\Delta$'s.
In other words, given $\Delta$ and a symmetric matrix, their flows
converge to a symmetric matrix with nonzero pattern $\Delta$.
We do not
know if the zero-preserving properties of these flows have been
studied.
Driessel~[2004] generalizes the Toda flow in a different way than we
do here.
\end{Remark}

We want to describe a geometrical explanation for the zero-preserving
property of the Toda flow. This geometrical reason apparently is not
well-known. (See, however, Symes~[1980a, 1980b, 1982].) Let $Upper(n)$
denote the group of invertible upper triangular matrices; in symbols,
$$
Upper(n) := \{ U \in Gl(n) : i>j \Rightarrow U(i,j) = 0 \}
$$
where $Gl(n)$ denotes the group of invertible $n\times n$ matrices.
Let $upper(n)$ denote the linear space of upper triangular matrices; in
symbols,
$$
upper(n) := \{ R\in \R^{n\times n}: i>j \Rightarrow R(i,j)=0 \}.
$$
Note that the space tangent to the matrix group of invertible upper
triangular matrices at the identity may be identified with the space
of upper triangular matrices; in symbols,
$$
Tan.Upper(n).I = upper(n) .
$$
Note that the space of square matrices $\R^{n \times n}$ is the direct
sum of the space of symmetric matrices and the space of strictly upper
triangular matrices since
$$
X = X_l + X_d + X_u = ( X_l + X_d + X_l^T) + ( X_u - X_l^T) .
$$
Let $\sigma: \R^{n \times n} \to Sym(n)$ denote the corresponding projection;
in symbols,
$$
\sigma.X := X_l + X_d + X_l^T .
$$
We consider the following map:
$$
\alpha := Upper(n) \times Sym(n) \to Sym(n) : (U,X) \to \sigma(UXU^{-1}) .
$$

\begin{Proposition}
The mapping $\alpha$ is a group action.
\end{Proposition}
\begin{Proof}
Note $\alpha(U_1,\alpha(U_2,X))=\alpha(U_1U_2,X)$ iff
$$
\sigma(U_1 \sigma(U_2XU_2^{-1}) U_1^{-1})
= \sigma(U_1 U_2 X U_2^{-1} U_1^{-1}).
$$
Let $Y$ be defined by $\sigma(U_2XU_2^{-1}) + Y := U_2XU_2^{-1}$.
Note that $Y$ is strictly upper-triangular. Then
$$
U_1(\sigma(U_2 X U_2^{-1}))U_1^{-1} + U_1 Y U_1^{-1}
= U_1U_2XU_2^{-1}U_1^{-1}.
$$
Note $\sigma(U_1 Y U_1^{-1}) = 0 $ since $U_1 Y U_1^{-1}$ is strictly
upper triangular.
\end{Proof}
For any symmetric matrix $X$, we have the orbit of $X$ under this
action:
$$
Orbit(X) = \alpha.Upper(n).X
   = \{ \sigma(UXU^{-1}) : U \in Upper(n)\} .
$$
Consider the following map determined by $X$:
$$
\beta.X := Upper(n) \to Sym(n): U \mapsto \sigma(UXU^{-1}) .
$$
Note that the image of this map is the orbit of $X$. We
differentiate $\beta.X$ to get the following linear map:
$$
D(\beta.X).I = Tan.Upper(n).I \to Tan.Sym(n).X : R \mapsto \sigma[R,X] .
$$
As noted above we can identify $Tan.Upper(n).I$ with $upper(n)$.
As before, we can also identify $Tan.Sym(n).X$ with $Sym(n)$. We shall
regard $D(\beta.X).I$ as a linear map from $upper(n)$ to $Sym(n)$. It is
not hard to prove that the space tangent to the orbit of $X$ at $X$
is the image of this linear map; in symbols,
$$
Tan.Orbit(X).X = \{\sigma[R,X] : R \in upper(n) \} .
$$
(For details see Warner[1983]
chapter 3: Lie groups, section: homogeneous manifolds.)

Let $T$ denote the tridiagonal symmetric matrix determined by the triple
$(1,0,1)$; in symbols,
$$
T:= \begin{pmatrix}
0 & 1 & 0 & \hdots & 0 & 0 \\
1 & 0 & 1 & \hdots & 0 & 0 \\
0 & 1 & 0 & \hdots & 0 & 0 \\
\vdots & \vdots & \vdots & \ddots & \vdots & \vdots \\
0 & 0 & 0 & \hdots & 0 & 1 \\
0 & 0 & 0 & \hdots & 1 & 0
\end{pmatrix} .
$$
We find the following result rather surprising. (In particular, we
do not know the historical origin of this result.)

\begin{Proposition}
The tridiagonal symmetric matrices with trace
equal zero and nonzero sub-diagonal (and super-diagonal) entries are
the orbit of the matrix $T$ under the action by the group $Upper(n)$.
\end{Proposition}

\flushpar
(The tridiagonal matrices with nonzero sub-diagonal and super-diagonal
are often called {\it Jacobi matrices}.)

\begin{Proof}
Note that every element of $Upper(n)$ can be written as the product of
an invertible diagonal matrix and an element of $Upper(n)$ with only
ones on the diagonal.
We sketch the rest of the proof when $n=3$; it should be clear how to
generalize these calculations. We use $*$ to denote irrelevant entries
in matrices. We have

\begin{align*}
&\begin{pmatrix}
1 & l_1 & \ast \\
0 & 1   & l_2  \\
0 & 0   &  1
\end{pmatrix}
\begin{pmatrix}
0 & 1 & 0 \\
1 & 0 & 1 \\
0 & 1 & 0
\end{pmatrix}
\begin{pmatrix}
1 & -l_1 & \ast \\
0 &  1   & -l_2 \\
0 &  0   &  1
\end{pmatrix}
=
\begin{pmatrix}
1 & l_1 & \ast \\
0 & 1   & l_2  \\
0 & 0   &  1
\end{pmatrix}
\begin{pmatrix}
0 & \ast & \ast \\
1 & -l_1 & \ast \\
0 &  1   & -l_2
\end{pmatrix}\\
=&
\begin{pmatrix}
l_1 & \ast    & \ast \\
1   & l_2-l_1 & \ast \\
0   &  1      & -l_2
\end{pmatrix}\\
\intertext{and}
&\begin{pmatrix}
d_1 & 0   & 0 \\
0   & d_2 & 0 \\
0   & 0   & d_3
\end{pmatrix}
\begin{pmatrix}
a_1 & 1      & 0 \\
1   & a_2    & 1 \\
0   & 1      & a_3
\end{pmatrix}
\begin{pmatrix}
1/d_1 & 0     &   0 \\
0     & 1/d_2 &   0 \\
0     & 0     & 1/d_3
\end{pmatrix} \\
=&
\begin{pmatrix}
a_1     &  \ast   & \ast \\
d_2/d_1 & a_2     & \ast \\
0       & d_3/d_2 & a_3
\end{pmatrix} .
\end{align*}
\end{Proof}

Note that for any symmetric matrix $X$, we have
$$
[X,X_l - X_l^T] = - [X, X_d + 2 X_l^T]
$$
since
\begin{align*}
0 &= [X,X] = [X, X_l + X_d + X_l^T] \\
  &= [X,X_l - X_l^T] + [X, X_d + 2 X_l^T] .
\end{align*}
Thus we see that the Toda flow can be rewritten as the following
``differential algebraic '' initial value problem:
\begin{align*}
X^\prime &= [X, X_l - X_l^T], \quad X(0) = A\\
[X,X_l-X_l^T] &= - [X, X_d + 2 X_l^T] .
\end{align*}
From the second equation (which holds trivially) we see that not only
does the solution stay on the iso-spectral surface, but it also stays
on the orbit of $A$ under the action by the upper triangular group.
Thus if $A$ is a tridiagonal matrix with trace zero and $X(t)$ is the
solution of the differential equation at time $t$ then $X(t)$ is
tridiagonal and has trace zero. It is easy to see that the zero trace
condition can be replaced by a constant trace condition.

It is not hard to see that these observations concerning symmetric
tridiagonal matrices generalize to any symmetric staircase pattern of
interest. The zero-preserving iso-spectral flow that
we derive in the main part of this report can be viewed as a differential 
algebraic equation
similar to the one we have here.

The Toda flow is also related to an optimization problem closely
related to the one we mentioned at the end of the introduction. As
there, let $D$ be a symmetric matrix and let
$f:=Sym(n)\to \R: X \mapsto (1/2)\langle X-D, X-D \rangle $ be
an ``objective function''. Consider the following optimization
problem:

\begin{Problem}
Given a symmetric matrix $A$, minimize $f(X)$ subject to the
constraint that $X$ is in $Iso(A)$.
\end{Problem}
This problem is analyzed in Chu and Driessel~[1990]. (See also
Driessel~[2004].)

Computing the derivative of $f$ we get that, for any symmetric matrices
$X$ and $H$, $DfX.H=\langle X-D,H \rangle$. For the gradient of $f$ at
$X$, we then have $\nabla f.X=X-D$. We can get an iso-spectral vector
field by orthogonal projection as follows. Let $l.X := D(\omega.X).I$.
Recall that, for all skew-symmetric $K$, $D(\omega.X).I.K=[K,X]$.  Note
that the adjoint $(l.X)^*$ of $l.X$ is the following map:
$$
(l.X)^* = Sym(n) \to Skew(n): Y \mapsto [Y,X]
$$
since, for every symmetric matrix $Y$ and every skew-symmetric
matrix $K$, $\langle [K,X],Y \rangle = \langle K,[Y,X] \rangle$.
If $l.X$ is injective then the projection onto $Tan.Iso(X).X$ is
the operator $(l.X) ((l.X)(l.X)^*)^{-1} (l.X)^*$.  Instead of
using this orthogonal projection, we simply use the map
$(l.X)(l.X)^*$; in other words, we drop the factor involving the
inverse from the projection formula. (For more on this matter see
Driessel~[2004].) We can also then drop the requirement that $l.X$
be injective. Note that $(l.X)(l.X)^*Y=[[Y,X],X]$. In particular,
we have $(l.X)(l.X)^*(-\nabla.f.X)= [[D-X,X],X] = [[D,X],X]$.  We
use a ``quasi-projection'' similar to this one in order to derive
our zero-preserving iso-spectral flow.

The {\bf double-bracket flow} is the flow associated with the
following differential equation:
$$
X^\prime = [[D,X],X] .
$$
\begin{Proposition} The double-bracket flow has the following properties:
\item{1.} This flow preserves eigenvalues.
\item{2.} The objective function $f$ is non-increasing along solutions of
this flow.
\item{3.} A symmetric matrix is an equilibrium point of this flow iff
it commutes with $D$.
\item{4.} Let $D$ be the diagonal matrix with diagonal entries
$1,2,\dots,n$. Then, on the space of tridiagonal symmetric matrices,
this flow coincides with the Toda flow.
\end{Proposition}

\begin{Proof} For any solution $X(t)$ of the double bracket differential
equation, we have
\begin{align*}
(f.X)^\prime &= \langle X-D,X^\prime \rangle
= \langle X-D,[[D,X],X] \rangle \\
&= \langle [X-D,X],[D,X] \rangle = - \langle [D,X],[D,X] \rangle \leq 0.
\end{align*}
This inequality shows the $f$ is non-increasing along solutions of this flow.
We leave the rest of the proof to the reader.
\end{Proof}


\section*{Acknowledgments.}
\addcontentsline{toc}{section}{Acknowledgments.}

We wrote most of this report during the Fall of 2001 while visiting the Fields
Institute for Research in Mathematical Sciences in Toronto, Ontario. We
wish to thank all the people at the institute who extended hospitality
to us during this pleasant visit, especially Ken Jackson (one of the
organizers of the Thematic Year on Numerical and Computational
Challenges in Science and Engineering at the Institute) and Ken
Davidson (head of the institute). In addition to these people, we wish
to thank colleagues who discussed iso-spectral flows with us during
this visit and offered advice and encouragement: Chandler Davis
(University of Toronto), Itamar Halevy (University of Toronto), Peter
Miegom (Fields Institute), and John Pryce (Cranfield University, UK).

Further, Alf Gerisch acknowledges financial support from the Fields
Institute for Research in Mathematical Sciences and the University of
Guelph.

In 2003, we received constructive comments from an anonymous
referee.  We thank the referee for these. In particular, the
suggestion that we compare projections with quasi-projections lead us to add
the appendix (based on our research notes of 2001) concerning this
matter.

\section*{References.}
\addcontentsline{toc}{section}{References.}

\flushpar Anderson, W. N., Jr. and Duffin, R. J. [1969] {\it Series and 
Parallel Addition of Matrices}, J. of Mathematical Analysis and Applications 
26, pp.~576-594.

\flushpar Anderson, W. N., Jr. [1971] {\it Shorted Operators}, 
SIAM J. Appl. Math. 20,~pp. 520-525.

\flushpar Anderson, W. N., Jr. and Schreiber, M. [1972] {\it The infimum of 
two projections}, Acta Sci. Math. 33, pp.~165-168. 

\flushpar Anderson, W. N., Jr. and Trapp, G. E. [1975] {\it Shorted Operators II},
SIAM J. Appl. Math. 28, pp.~60-71.

\flushpar Arbenz, P. and Golub, G. [1995] {\it Matrix shapes invariant
under the symmetric QR algorithm}, Numerical Lin. Alg. with
Applications 2, pp.~87-93.

\flushpar
Ashlock, D. A.; Driessel, K. R. and Hentzel, I. R. [1997a]
{\it Matrix structures invariant under Toda-like iso-spectral flows},
Lin. Alg. and Applications 254, pp.~29-48.

\flushpar
Ashlock, D. A.; Driessel, K. R. and Hentzel, I. R. [1997b]
{\it Matrix structures invariant under Toda-like iso-spectral flows:
sign-scaled algebras}, preprint.

\flushpar Bellman, R. [1970] Introduction to Matrix Analysis,
McGraw-Hill.

\flushpar
Chu, M. and Driessel, K. R. [1990] {\it The projected gradient method for
least squares matrix approximations with spectral constraints}, SIAM
J. Numer. Anal. 27, pp.~1050-1060.

\flushpar
Chu, M. and Norris, L.K. [1988], {\it Iso-spectral flows and abstract
matrix factorizations}, SIAM J. Numer. Anal. 25, pp.~1383-1391.

\flushpar
Curtis, M.L. [1984], Matrix Groups, Springer-Verlag.

\flushpar
Demmel, J. [1997] Applied Linear Algebra, SIAM, Section 5.5: Differential
Equations and Eigenvalue Problems.

\flushpar
Driessel, K.R. [2004] {\it Computing canonical forms using flows},
Linear Algebra and its Applications 379, pp.~353-379.

\flushpar Fasino, D.[2001], {\it Iso-spectral flows on
displacement structured matrix spaces}, in Structured Matrices:
Recent Developments in Theory and Computation, D. Bini, E
Tyrtyshnikov and P. Yalamov (editors), Nova Science Publisher Inc.

\flushpar
Halmos, P.R. [1958], Finite-Dimensional Vector Spaces, D. Van
Nostrand, Inc.

\flushpar
Hirsch, M.W.; Smale, S. and Devaney, R.L. [2004],
Differential Equations, Dynamical Systems \& An Introduction to Chaos,
2nd Edition, Elsevier.

\flushpar
Kubo, K. and Ando, T. [1980] {\it Means of positive linear operators},
Mathematische Annalen 246, pp.~205-224.

\flushpar
Lawson, C. and Hanson, R. [1974], Solving Least Squares Problems,
Prentice-Hall, Inc.

\flushpar
Leon, S. [1986], Linear Algebra with Applications, Macmillan
Publishing Company.

\flushpar
Strang, G. [1980], Linear Algebra and Its Applications, Academic Press.

\flushpar
Symes, W. W. [1980a] {\it Systems of Toda type, inverse spectral
problems, and representation theory}, Inventiones Math. 59, pp.~13-51.

\flushpar
Symes, W. W. [1980b] {\it Hamiltonian group actions and integrable
systems}, Physica 1D, pp.~339-374.

\flushpar
Symes, W. W. [1982] {\it The QR algorithm and scattering for the finite
non-periodic Toda lattice}, Physica 4D, pp.~275-280.

\flushpar
Warner, F.W. [1983] Foundations of Differential
Manifolds and Lie Groups, Springer-Verlag.

\end{document}